\documentclass{amsart}

\usepackage{amsthm}
\usepackage{amssymb}
\usepackage{amsmath}
\newtheorem {Theorem}   {Theorem} 

\numberwithin{Theorem}{section}

\theoremstyle{definition}

\theoremstyle{remark}

\newtheorem {Corollary}[Theorem]{Corollary}

\begin{document}

\title{Flows on solenoids are generically not almost periodic}

\author{Alex Clark}

\address{Department of Mathematics, University of North Texas, 
Denton, TX 76203-5118}
\email{alexc@unt.edu}

\date{January 10, 1999 and, in revised form, June 20, 1999.}

\thanks{This work was funded in part by a faculty
research grant from the University of North Texas.}


\begin{abstract}

The space of non--singular flows on the solenoid $\Sigma _{N}$ is shown to
contain a dense $G_{\delta }$ consisting of flows which are not almost
periodic. Whether this result carries over to
Hamiltonian flows remains an open question.

\end{abstract}

\maketitle

\section*{Introduction}

For any compact symplectic manifold $M$ of dimension at least $4$, L. Markus
and K. R. Meyer demonstrate in \cite{MM} that the space $\mathfrak{H}^{k}\left(
M\right) $ of all $C^{k}$ ($k\geq 4$) Hamiltonians on $M$ contains a generic
subset $\mathfrak{M}_{\Sigma }$ such that for each Hamiltonian $dH^{\#}$ in $
\mathfrak{M}_{\Sigma }$ and for each solenoid $\Sigma _{N}$ there exists a
minimal set for the flow induced by $dH^{\#}$ that is homeomorphic to $
\Sigma _{N}$. They leave open the question of whether the flows on these
solenoids are almost periodic (\cite{MM}, p. 90).

Corresponding to a sequence of natural numbers $N=\left(
n_{1},n_{2},...\right) $ with $n_{j}\geq 2$ for each $j$ there is the
solenoid $\Sigma _{N}$ which is\textbf{\ }the inverse limit of the inverse
sequence $\left\{ X_{j},f_{j}\right\} _{j=1}^{\infty }$ with factor space $
X_{j}=S^{1}=\{\exp \left( 2\pi it\right) \in \mathbb{C\mid }t\in \lbrack 0,1)\}$
for each $j$ and bonding maps $f_{j}\left( z\right) =z^{n_{j}}$
\[
S^{1}\stackrel{n_{1}}{\longleftarrow }S^{1}\stackrel{n_{2}}{\longleftarrow }
S^{1}\stackrel{n_{3}}{\longleftarrow }\cdots \hspace{0.2in}\Sigma
_{N}=\left\{ \left\langle z_{j}\right\rangle _{j=1}^{\infty }\in
\prod_{i=1}^{\infty }S^{1}\mid z_{j}=f_{j}\left( z_{j+1}\right) \text{ for }
j=1,2,...\right\} . 
\]
Here $\prod_{i=1}^{\infty }S^{1}$ is the compact topological group with the
group operation (written $``+"$) given by factor--wise multiplication and
with the metric 
\[
d\left( \left\langle x_{j}\right\rangle _{j=1}^{\infty },\left\langle
y_{j}\right\rangle _{j=1}^{\infty }\right) =\sum_{j=1}^{\infty }\frac{1}{
2^{j}}\left| x_{j}-y_{j}\right| . 
\]
This metric\ also serves as a metric for the subgroup $\Sigma _{N}$. There
is then the continuous (but not bicontinuous) isomorphism $\pi _{N}:\mathbb{
R\rightarrow }C_{N}$ onto the $\Sigma _{N}$--arc component $C_{N}$ of the
identity $e=\left( 1,1,...\right) $ given by
\[ 
t\mapsto \left\langle \exp
\left( 2\pi it\right) ,\exp \left( \dfrac{2\pi it}{n_{1}}\right) ,...,\exp
\left( \dfrac{2\pi it}{n_{1}\cdots n_{j}}\right) ,...\right\rangle ,
\] 
and we have the family of linear flows $\Lambda _{N}$ on $\Sigma _{N}$: 
\begin{eqnarray*}
\Lambda _{N} &=&\left\{ \phi _{N}^{\alpha }:\mathbb{R}\times \Sigma
_{N}\rightarrow \Sigma _{N}\mathbb{\mid \alpha \in R}\right\} \\
\phi _{N}^{\alpha }\left( t,x\right) &=&\pi _{N}\left( \alpha t\right) +x
\end{eqnarray*}
and any almost periodic flow on $\Sigma _{N}$ is equivalent ($C^{0}$
conjugate) to some $\phi _{N}^{\alpha }\in \Lambda _{N}$ \cite{C}. (A \emph{flow}
is a continuous group action of $\left( \mathbb{R},+\right) $.)

Given any point $x\in \Sigma _{N}$ the map $\pi _{N}+x$ sending $t\mapsto
\pi _{N}\left( t\right) +x$ parameterizes the arc component of $x$, which
will coincide with the trajectory of $x$ for any non--singular flow $\phi $
on $\Sigma _{N}$. For $k=0,1,...,$ $\infty $ and $0\leq \alpha \leq 1$, we
consider a flow $\phi :\mathbb{R}\times \Sigma _{N}\rightarrow \Sigma _{N}$ to
be $C^{k+1+\alpha }$ if the associated ``vector field'' 
\[
v_{\phi }:\Sigma _{N}\mathbb{\rightarrow R};\;x\longmapsto \frac{d\left[ \left(
\pi _{N}+x\right) ^{-1}\phi \left( t,x\right) \right] }{dt}\mid _{t=0} 
\]
exists and satisfies the conditions:

\begin{enumerate}
\item  for $j=1,...,k$ the function $v_{\phi }^{j}:\Sigma _{N}\mathbb{
\rightarrow R};\;x\longmapsto \dfrac{d^{j}\left[ \left( \pi _{N}+x\right)
^{-1}\phi \left( t,x\right) \right] }{dt^{j}}\mid _{t=0}$ is continuous [$
v_{\phi }^{1}=v_{\phi }$]

\item  for each $x\in \Sigma _{N}$ the function $\mathbb{R\rightarrow R}
;\;t\longmapsto v_{\phi }^{k}\left( x+\pi _{N}\left( t\right) \right) $ is $
\alpha $--H\"{o}lder (see below) if $\alpha >0$ and $k<\infty $.
\end{enumerate}

\noindent We then endow the space $C^{k+\alpha }\left( N\right) $ of vector
fields $v_{\phi }$\ on $\Sigma _{N}$ stemming from $C^{k+1+\alpha }$ flows $
\phi $ with the metric 
\[
d_{k}\left( v_{\phi },v_{\psi }\right) \stackrel{def}{=}\sum_{j=1}^{k}
\max_{x\in \Sigma _{N}}\left| v_{\phi }^{j}\left( x\right) -v_{\psi
}^{j}\left( x\right) \right| \text{ for }1\leq k<\infty 
\]
\[
\text{and }d_{\infty }\left( v_{\phi },v_{\psi }\right) \stackrel{def}{=}
\sum_{r=1}^{\infty }\frac{2^{-r}d_{r}\left( v_{\phi },v_{\psi }\right) }{
1+d_{r}\left( v_{\phi },v_{\psi }\right) }. 
\]
If $\Sigma _{N}$ is embedded in some compact $C^{k}$ manifold so that the
flow $\phi _{N}^{1}$ extends to a $C^{k}$ flow on the manifold, then this
metric induces the same topology as the restriction of the Whitney topology
(see \cite{MM}, pp. 34--35 for a description of the Whitney topology). We are
interested in the non--singular flows on $\Sigma _{N}$ and so will work with 
$C^{k+\alpha }\left( N+\right) ,$ the subspace of $C^{k+\alpha }\left(
N\right) $ consisting of positive $v_{\phi }$. (Whenever $\phi $ is
non--singular the map $v_{\phi }$ must be either positive or negative since
each trajectory is dense and the intermediate value theorem yields a
singularity if there is a change of signs, and so we consider positive $
v_{\phi }$ without loss of generality.)

If $p_{j}:\Sigma _{N}\rightarrow S^{1}$ denotes the projection onto the $
j^{th}$ factor $\left\langle x_{k}\right\rangle _{k=1}^{\infty }\mapsto
x_{j} $ and if $\widehat{N}$\ denotes the subgroup of $\left( \mathbb{R}
,+\right) $ generated by 
\[
\left\{ \dfrac{1}{n_{1}},\dfrac{1}{n_{1}n_{2}},...,
\dfrac{1}{n_{1}n_{2}\cdots n_{j}}\mid j=1,2,...\right\} 
\] 
and if for $r=
\dfrac{m}{n_{1}n_{2}\cdots n_{j}}\in \widehat{N}$\ \ $\chi _{r}:\Sigma
_{N}\rightarrow S^{1}$ denotes the homomorphism sending $x\mapsto \left(
p_{j+1}\left( x\right) \right) ^{m},$ then $\Xi \left( N\right) =\left\{ \chi
_{r}\mid r\in \widehat{N}\right\} $ is the group of characters of $\Sigma
_{N}$. Given a continuous function $f:\Sigma _{N}\mathbb{\rightarrow C}$, the
theorem of Peter and Weyl guarantees that $f$ can be uniformly approximated
by finite linear combinations of characters $\chi _{r}$ (see, e.g., \cite{P}), and
since $\chi _{\frac{m}{n}}\circ \pi _{N}\left( t\right) =\exp \left( 2\pi
imt/n\right) $ the map $f_{e}:\mathbb{R\rightarrow C};$ $t\mapsto f\left( \pi
_{N}\left( t\right) \right) $ is the uniform limit of periodic maps and is
thus\emph{\ limit periodic} by definition (see, e.g., \cite{Be} 1\S 6), which is a
special type of almost periodic function. It then follows that the
Bohr--Fourier series 
\[
f_{e}\left( t\right) \sim \sum_{r\in \hat{N}}f_{r}\chi _{r}\circ \pi
_{N}\left( t\right) =\sum_{r\in \hat{N}}f_{r}\exp \left( 2\pi irt\right) 
\]
when appropriately ordered converges uniformly to $f_{e}$ provided $f_{e}$
is $\alpha $--H\"{o}lder for some $0<\alpha \leq 1$: 
\[
\sup_{t}\left| f_{e}\left( t+\delta \right) -f_{e}\left( t\right) \right|
<C\delta ^{\alpha }\text{\ }\left( \text{for some }C>0\text{ and all }\delta
>0\right) \text{ (see \cite{Be} 1\S 8).} 
\]
Whenever we are dealing with such a function we will assume without comment
that the series is so arranged that $f=\sum\limits_{r\in \hat{N}}f_{r}\chi
_{r}$.

We shall show that in each space $C^{k+\alpha }\left( N+\right) $ the
collection of vector fields $v_{\phi }$ corresponding to flows $\phi $ which
are not almost periodic contains a dense $G_{\delta }$ provided $k\geq 1$ or 
$\alpha =1$ (in other words: $v_{\phi }$ is Lipschitz).

\section{Proof of the Main Result}

Fix a space $C^{k+\alpha }\left( N+\right) $ with $k\geq 1$ or $\alpha =1$.\
As already mentioned, a flow $\phi $ on $\Sigma _{N}$ is almost periodic
exactly when there is a homeomorphism $h^{\prime }:\Sigma _{N}\rightarrow
\Sigma _{N}$ providing an equivalence of $\phi $ with some linear flow $\phi
_{N}^{\alpha ^{\prime }}$, $h^{\prime }\circ \phi \left( t,x\right) \equiv
\phi _{N}^{\alpha ^{\prime }}\left( t,h^{\prime }\left( x\right) \right) $.
Any such $h^{\prime }$ is homotopic to an automorphism $a$ of $\Sigma _{N}$,
and the automorphism $a^{-1}$ will in turn equate $\phi _{N}^{\alpha
^{\prime }}$ with a linear flow $\phi _{N}^{\alpha }$ (see \cite{S} and \cite{C}). And so 
$\phi $ is almost periodic exactly when there is a homeomorphism $
h=a^{-1}\circ h^{\prime }$ homotopic to the identity providing an
equivalence between $\phi $ and a linear flow $\phi _{N}^{\alpha }$.

We now translate the existence of such an $h$ into a more useful form, so
assume $\phi $ is almost periodic and that such an $h$ exists. Each $x\in
\Sigma _{N}$ belongs to the section $\mathcal{S}_{x}\stackrel{def}{=}
p_{1}^{-1}\left( p_{1}\left( x\right) \right) $ (which is a Cantor set), and
the time it takes the linear flow $\phi _{N}^{\alpha }$ to return to $
\mathcal{S}_{x}$ is the constant $\dfrac{1}{\alpha }$ since 
\[
p_{1}\left( \phi _{N}^{\alpha }\left( t,x\right) \right) =p_{1}\left( \pi
_{N}\left( t\alpha \right) +x\right) =p_{1}\left( \pi _{N}\left( t\alpha
\right) \right) \cdot p_{1}\left( x\right) =\exp \left( 2\pi it\alpha
\right) \cdot p_{1}\left( x\right) 
\]
and the smallest positive value of $t$ satisfying $\exp \left( 2\pi it\alpha
\right) =1$\ is $\dfrac{1}{\alpha }$. Since $1/x$ is $C^{\infty }$ on any
closed interval not containing $0$, $\lambda \left( x\right) \stackrel{def}{=
}\dfrac{1}{v_{\phi }\left( x\right) }=\sum\limits_{r\in \hat{N}}\lambda
_{r}\chi _{r}\left( x\right) $ is $C^{k+\alpha }$ and the $\phi $--return
time $\tau $ of $x$ to the section $\mathcal{S}_{x}$ is given by 
\begin{eqnarray*}
\tau \left( x\right) &=&\int_{0}^{1}\lambda \left( x+\pi _{N}\left( t\right)
\right) dt\\
&=&\lambda _{0}+\sum\limits_{r\in \hat{N}-\left\{ 0\right\} }\frac{
\lambda _{r}}{2\pi ir}\chi _{r}\left( x+\pi _{N}\left( 1\right) \right)
-\sum\limits_{r\in \hat{N}-\left\{ 0\right\} }\frac{\lambda _{r}}{2\pi ir}
\chi _{r}\left( x\right) \\
&\sim &\lambda _{0}+\sum\limits_{r\in \hat{N}-\left\{ 0\right\} }\frac{
\lambda _{r}}{2\pi ir}\left[ \chi _{r}\left( \pi _{N}\left( 1\right) \right)
-1\right] \chi _{r}\left( x\right) =\sum_{r\in \hat{N}}\tau _{r}\chi
_{r}\left( x\right) .
\end{eqnarray*}
Now $h$ is homotopic to the identity and so we have the map $\delta \sim
\sum\limits_{r\in \hat{N}}\delta _{r}\chi _{r}:\Sigma _{N}\rightarrow \mathbb{R}
$
\[
\delta \left( x\right) \stackrel{def}{=}\text{the unique time }t\text{
satisfying }h\left( x\right) =\phi _{N}^{\alpha }\left( t,x\right) 
\]
since $h\left( x\right) $ must lie in the same path component as $x$. Since $
h$ provides a flow equivalence, the $\phi _{N}^{\alpha }$--return time to
the $\phi _{N}^{\alpha }$--section $h\left( \mathcal{S}_{x}\right) $ for $
h\left( x\right) $ is $\tau \left( x\right) $. We can now express the
constancy of the $\phi _{N}^{\alpha }$--return time to $\mathcal{S}_{x}$ as 
\[
\delta \left( x\right) +\tau \left( x\right) -\delta \left( x+\pi _{N}\left(
1\right) \right) =\frac{1}{\alpha } 
\]
or 
\[
\tau \left( x\right) -\frac{1}{\alpha }=\delta \left( x+\pi _{N}\left(
1\right) \right) -\delta \left( x\right) . 
\]
We also have 
\[
\delta \left( x+\pi _{N}\left( 1\right) \right) -\delta \left( x\right) \sim
\sum\limits_{r\in \hat{N}-\left\{ 0\right\} }\left[ \chi _{r}\left( \pi
_{N}\left( 1\right) \right) -1\right] \delta _{r}\chi _{r}\left( x\right) 
\]
leading to the conditions 
\[
\delta _{r}=\frac{\tau _{r}}{\left[ \chi _{r}\left( \pi _{N}\left( 1\right)
\right) -1\right] }=\frac{\lambda _{r}}{2\pi ir}\text{ for }r\in \hat{N}
-\left\{ 0\right\} 
\]
and 
\[
\lambda _{0}=\tau _{0}=\frac{1}{\alpha }. 
\]
(There is no restriction on $\delta _{0}$ since we may follow $h$ by
translations of elements in $C_{N}$ and obtain other flow equivalences
homotopic to the identity.)

We then have the limit periodic function 
\[
\delta _{e}\left( t\right) \sim \sum\limits_{r\in \hat{N}}\delta _{r}\exp
\left( 2\pi irt\right) =\delta _{0}+\sum\limits_{r\in \hat{N}-\left\{
0\right\} }\frac{\lambda _{r}}{2\pi ir}\exp \left( 2\pi irt\right) . 
\]
According to Bohr's theorem on the integral of an almost periodic function
(see \cite{B} \S 68--69), the integral $F\left( T\right) =\int_{0}^{T}f\left(
t\right) dt$ of a limit periodic function $f\left( t\right) \sim
\sum\limits_{r\in \hat{N}}f_{r}\exp \left( 2\pi irt\right) $ is almost
periodic if and only if $F$ is bounded, in which case 
\[
F\left( t\right) \sim \sum\limits_{r\in \hat{N}}F_{r}\exp \left( 2\pi
irt\right) =F_{0}+\sum\limits_{r\in \hat{N}-\left\{ 0\right\} }\frac{f_{r}}{
2\pi ir}\exp \left( 2\pi irt\right) . 
\]
Comparing the Fourier--Bohr series of $\delta _{e}$ with that of $\lambda
_{e}$, we see that the bounded function $\delta _{e}\left( t\right) $
represents an integral of $\lambda _{e}\left( t\right) -\lambda _{0}$.
Moreover, whenever $\lambda _{e}\left( t\right) -\lambda _{0}$ has a bounded
integral we are able to construct a map $\delta $ as above, which in turn
allows us to construct an equivalence $h$. Notice that by its construction $
\delta $ will be as smooth as $\lambda $. This gives us the following result.

\begin{Theorem}
\bigskip The flow $\phi \in C^{k+\alpha }\left( N+\right) $ with vector
field $v_{\phi }$ is almost periodic if and only if $\lambda _{e}\left(
t\right) -\lambda _{0}$ has a bounded integral, where $\lambda _{e}\left(
t\right) =\dfrac{1}{v_{\phi }\circ \pi _{N}\left( t\right) }
=\sum\limits_{r\in \hat{N}}\lambda _{r}\exp \left( 2\pi irt\right) $.
\end{Theorem}

If $f\left( t\right) $ is periodic with Fourier series $\sum\limits_{r\in 
\hat{N}}f_{r}\exp \left( 2\pi irt\right) $, a necessary and sufficient
condition that it have a periodic integral is that 
\[
f_{0}=M\left\{ f\right\} 
\stackrel{def}{=}\lim\limits_{T\rightarrow \infty }\dfrac{1}{T}
\int_{0}^{T}f\left( t\right) dt=0 \text{ (see \cite{B} \S 68),}
\] 
and since any function $
\lambda _{e}\left( t\right) -\lambda _{0}$ as above can be approximated
(relative to $d_{k}$) arbitrarily closely by a periodic ``partial series''
(possibly containing infinitely many terms) of its Fourier--Bohr series and
since $M\left\{ f\right\} =0$ for any such partial series $f$,\ we obtain
the following.

\begin{Corollary}
The collection of almost periodic flows in $C^{k+\alpha }\left( N+\right) $
is dense.
\end{Corollary}

(We are making use of the fact that the function $v_{\phi }\mapsto 1/v_{\phi
}=\lambda $ is a homeomorphism of $C^{k+\alpha }\left( N+\right) .$) We
postpone the existence of flows for given functions until the next section.

Now we examine the functions $\lambda =1/v_{\phi }$ corresponding to flows $
\phi $ which are not almost periodic. First, we note 
\[
\left| M\left\{ f\right\} -M\left\{ g\right\} \right| \leq
\lim\limits_{T\rightarrow \infty }\dfrac{1}{T}\int_{0}^{T}\left| f\left(
t\right) -g\left( t\right) \right| dt\leq \sup \left| f\left( t\right)
-g\left( t\right) \right| . 
\]
And so $d_{1}\left( \lambda -\lambda _{0},\mu -\mu _{0}\right) \leq
2d_{1}\left( \lambda ,\mu \right) $. Now for $n=1,2,...$ we define the sets 
\[
U_{n}\stackrel{def}{=}\left\{ \lambda \mid \text{ there is a }T_{n}\text{
with }\left| \int_{0}^{T_{n}}\left( \lambda _{e}\left( t\right) -\lambda
_{0}\right) dt\right| >n\right\} 
\]
and claim that $U_{n}$ is open for $n=1,2,...$ . So suppose then that $
\lambda \in U_{n}$ with 
\[
\delta =\left| \int_{0}^{T_{n}}\left( \lambda _{e}\left( t\right) -\lambda
_{0}\right) dt\right| -n>0. 
\]
Now if $d_{1}\left( \lambda ,\mu \right) \leq d_{k}\left( \lambda ,\mu
\right) <\dfrac{\delta }{3\left| T_{n}\right| }$ we have 
\[
\left| \int_{0}^{T_{n}}\left( \lambda _{e}\left( t\right) -\lambda
_{0}\right) dt-\int_{0}^{T_{n}}\left( \mu _{e}\left( t\right) -\mu
_{0}\right) dt\right| \leq \left| \int_{0}^{T_{n}}\dfrac{2\delta }{3\left|
T_{n}\right| }dt\right| <\delta 
\]
and so $\mu $ is also in $U_{n}$, demonstrating that $U_{n}$ is open. And so
the $G_{\delta }$ $\cap _{n=1}^{\infty }U_{n}$ is the collection of $\lambda
=1/v_{\phi }$ corresponding to flows $\phi $ which are not almost periodic.

It then remains to show that $\cap _{n=1}^{\infty }U_{n}$ is dense. Let $
\lambda =1/v_{\phi }$ be given. We need to find an arbitrarily close
function which corresponds to a flow which is not almost periodic. Consider
for $m=1,2,...$ the $C^{\infty }$ maps $\Sigma _{N}\rightarrow \mathbb{R}$
\[
\rho _{m}\left( x\right) \stackrel{def}{=}\sum_{j=m}^{\infty }\frac{1}{
n_{1}\cdots n_{j}}\chi _{\frac{1}{n_{1}\cdots n_{j}}}\left( x\right) . 
\]
Bohr's theorem shows that $\rho _{m}\circ \pi _{N}\left( t\right) $ has an
unbounded integral since Parseval's equation for almost periodic functions
would fail for the Fourier--Bohr series of a bounded $\int_{0}^{T}\rho
_{m}\circ \pi _{N}\left( t\right) dt$, implying that at least one of $\text{
Re}\left( \rho _{m}\circ \pi _{N}\left( t\right) \right) $ and $\text{Im}
\left( \rho _{m}\circ \pi _{N}\left( t\right) \right) $ too has an unbounded
integral. (Notice that both the real and imaginary parts of a function with
mean value $0$ also have mean value $0$.) If $\lambda _{e}\left( t\right)
-\lambda _{0}$ has unbounded integral, there is nothing to prove; and if
not, we can choose $\lambda +\text{Re}\rho _{m}$ and $\lambda +\text{Im}\rho
_{m}$ positive and as close to $\lambda $\ as desired by choosing $m$ large
enough since $\left| \rho _{m}\left( x\right) \right| \leq \dfrac{1}{2^{m-1}}
$ (and similarly for the derivatives), and at least one of these two functions 
(say $\lambda ^{m}$) will be such that $\lambda ^{m}\circ \pi _{N}$ will have 
an unbounded integral when its mean value $\lambda _{0}$ is subtracted.

\begin{Corollary}
The collection of flows in $C^{k+\alpha }\left( N+\right) $ which are not
almost periodic is a dense $G_{\delta }$.
\end{Corollary}

\section{Realization of a flow for a given vector field}

Let $v$ be a positive Lipschitz function $\Sigma _{N}\rightarrow \mathbb{R}$.
We seek a flow $\phi _{v}$\ on $\Sigma _{N}$ which has $v$ as its vector
field. For $n\geq 2$ there is a $C^{\infty }$ flow $\phi $ on the tube $\Pi
_{0}=B^{2n-1}\times S^{1}$ ($B^{2n-1}$ is the closed unit ball of $\mathbb{R}
^{2n-1}$) satisfying:

\begin{enumerate}
\item  $\phi $ has a homeomorphic copy $\mathfrak{S}_{N}$ of $\Sigma _{N}$ as a
limit set

\item  $d\psi /dt\equiv 1$ along all the orbits of $\phi ,$ where $0\leq
\psi <1$ parameterizes the $S^{1}$ factor of $\Pi _{0}$ (see \cite{MM}, p. 87).
\end{enumerate}

The limit set $\mathfrak{S}_{N}$ is realized as the nested intersection of
compact tubes $\Pi _{j}$, $j=0,1,2,...$ , \ satisfying for $j=0,1,...$ :

\begin{enumerate}
\item  $\Pi _{j}$ is homeomorphic to $\Pi _{0}$ and

\item  $\Pi _{j+1}$ encircles $\Pi _{j}$\hspace{0.04in}$n_{j+1}$ times.
\end{enumerate}

Therefore, $\phi $ restricted to $\mathfrak{S}_{N}$ has the constant return time
of $1$ to the sections 
\[
\mathcal{S}_{t}\stackrel{def}{=}\left\{ \left( x,\psi \right) \in
B^{2n-1}\times S^{1}\cap \mathfrak{S}_{N}\mid \psi =t\right\} 
\]
and is thus equivalent to the linear flow $\phi _{N}^{1}$ on $\Sigma _{N}$
via a homeomorphism, say $h$, where without loss of generality $\psi(h^{-1}(e))=0 $. 
Associated with the flow $\phi $ is a vector
field $\mathfrak{v}$. We shall obtain the desired flow $\phi _{v}$ by a time
change of the flow $\phi $: we shall multiply the vectors of $\mathfrak{v}$ by a
function to change their lengths but not their directions and then use $h$
to obtain $\phi _{v}$. (Here we are measuring lengths of vectors by the lengths
of their projections onto the $\mathbb{R}$ factor in the tangent bundle corresponding 
to the $S^{1}$ factor in $\Pi _{0}$.) 

Setting $\psi_{0}=\psi$, for $j=1,2,... $ let $0\leq \psi _{j}<n_{1}\cdots n_{j}$
parameterize the $S^{1}$ factor of $\Pi _{j}$ in such a way that for 
$\left( x,\psi_{j}\right) =\left( x,\psi _{j-1}\right) $ in $\Pi _{j}$ 
we have $\psi_{j}=\psi _{j-1}\left( \text{mod}1\right) $ and $\psi_{j}(h^{-1}(e))=0 $.
(We could use $\phi $--sections transverse to the elliptic
periodic orbits $\gamma _{0},\gamma _{1},...$ for example (see \cite{MM}, p. 87).)
For $j=0,1,2,...$ let $s_{j}\left( t\right) $ be a purely periodic Lipschitz
partial sum of $v_{e}\left( t\right) =\lim\limits_{j\rightarrow \infty
}s_{j}\left( t\right) =\sum\limits_{r\in \hat{N}}v_{r}\exp \left( 2\pi
irt\right) $, where the period of $s_{0}$ is $1$ and the period of $s_{j}$
is $n_{1}\cdots n_{j}$ for $j=1,2,...$ (see \cite{Be} 1\S 8). We then form for $
j=0,1,2,...$ the following functions $\tau _{j}:\Pi _{0}\rightarrow \mathbb{R}$; $
\tau _{0}\left( \left( x,\psi _{0}\right) \right) =s_{0}\left( \psi
_{0}\right) $,\ $\tau _{1}\left( \left( x,\psi _{0}\right) \right) =\tau
_{0}\left( \left( x,\psi _{0}\right) \right) $ for $\left( x,\psi
_{0}\right) $ in $\Pi _{0}-\mathcal{T}_{1}$ (where $\mathcal{T}_{1}$ is a
tube satisfying $\Pi _{1}\subset \mathcal{T}_{1}\subset \Pi _{0}$) and $\tau
_{1}\left( \left( x,\psi _{0}\right) \right) $ is $\tau _{0}\left( \left(
x,\psi _{0}\right) \right) $ gradually changed in $\mathcal{T}_{1}-\Pi _{1}$
until finally $\tau _{1}\left( \left( x,\psi _{1}\right) \right)
=s_{1}\left( \psi _{1}\right) $ for $\left( x,\psi _{1}\right) \in \Pi _{1}$
. Continue in the same manner so that $\tau _{j+1}\left( \left( x,\psi
_{0}\right) \right) =\tau _{j}\left( \left( x,\psi _{0}\right) \right) $ for 
$\left( x,\psi _{0}\right) $ in $\Pi _{0}-\mathcal{T}_{j+1}$ (where $
\mathcal{T}_{j+1}$ is a tube satisfying $\Pi _{j+1}\subset \mathcal{T}
_{j+1}\subset \Pi _{j}$) and $\tau _{j+1}\left( \left( x,\psi _{0}\right)
\right) $ is $\tau _{j}\left( \left( x,\psi _{0}\right) \right) $ gradually
changed in $\mathcal{T}_{j+1}-\Pi _{j+1}$ until finally $\tau _{j+1}\left(
\left( x,\psi _{j+1}\right) \right) =s_{j+1}\left( \psi _{j+1}\right) $ for $
\left( x,\psi _{j+1}\right) \in \Pi _{j+1}$. And then with $\tau \left(
\left( x,\psi _{0}\right) \right) \stackrel{def}{=}\lim\limits_{j\rightarrow
\infty }\tau _{j}\left( \left( x,\psi _{0}\right) \right) $ we have the
function to obtain the desired time change of $\mathfrak{v}$.

This also shows that for arbitrarily small time changes we can alter a
generic class of Hamiltonian flows to obtain flows which are not Lyapunov
stable on solenoidal minimal sets (see \cite{NS} V\S 8). By adjusting the coefficients of the $
\rho _{m}$ as in the previous section, we even have flows which have points
on the same orbit arbitrarily close that ``lap'' one another relative to
these tubes. However, it is not clear what happens if we restrict ourselves
to Hamiltonian flows. To make matters even more complicated, it is not clear
that these solenoids persist as limit sets under small perturbations (see \cite{AM}
8\S 5).

\end{document}